# Vortex Filament Equation and Non-Lightlike Space Curves in Minkowski 3-space

Nevin Gürbüz[1] and Süleyman Cengiz[2*]

[1]Department of Mathematics, Eskişehir Osmangazi University, Turkey

[2]Department of Mathematics, Çankırı Karatekin University, Turkey

E-mail: suleymancengiz@karatekin.edu.tr

**Abstract.** In this paper, the parallel transport frames over non-lightlike curves in Minkowski 3-space are introduced. Evolution equations of these frames with respect to arc length and time are calculated over the space of these curves. Then the equivalence of the non-linear Schrödinger equation and non-linear heat system to the vortex filament equation for the binormal motion of non-lightlike curves is shown by use of the parallel frames and their evolution equations.

**Keywords.** Vortex filament equation, non-lightlike curves, Minkowski 3-space, parallel frame, NLS, heat system

## 1. Introduction

The subject of vortex dynamics has been initiated by the seminal paper of Helmholtz in 1858. In the beginning of 20th century da Rios modeled the evolution of a thin vortex filament in a viscous fluid by the motion of a curve propagating in $R^3$ according to the vortex filament equation. The works of Da Rios and his master Levi-Civita on vortex filament motion and asymptotic potential theory remained unnoticed for a long time. At 1960's their work has been re-discovered by many researchers. In 1971, Hasimoto showed the equivalence of the vortex filament equation with the non-linear Schrödinger equation (NLS) (1972).

The intrinsic geometry of the nonlinear Schrödinger equation in Euclidean 3-space has been described by Rogers and Chief (1998).Then a similar work for the non-linear Schrödinger equation of repulsive type (NLS⁻) and the non-linear heat system was performed in a general intrinsic geometric setting including a normal congruence in Minkowski 3-space by Gürbüz (2010). The equivalence of the binormal motions of spacelike and timelike curves in Minkowski 3-space to the non-linear heat system and the NLS¯ equation was shown in (Ding, 1998; Ding and Inoguchi, 2004) by means of the Schrödinger flows into the spaces $S^{1,1}$ and $H^2$, respectively. The equivalence of these systems and the motions of curves in $R^3$ and $S^3$ was demonstrated by use of the parallel frames in (Terng, 2005; Zhang and Wu, 2007). In this paper, we extend this demonstration to Minkowski 3-space with use of the same idea as in (Terng, 2005).

---

[*] Corresponding Author

## 2. Preliminaries

The Minkowski 3-space $E_1^3$ is the Euclidean 3-space $E^3$ provided with the Lorentzian inner product

$$<x,y>_L = -x_1 y_1 + x_2 y_2 + x_3 y_3$$

where $x = (x_1, x_2, x_3), y = (y_1, y_2, y_3) \in E^3$. A vector $x \in E_1^3$ is spacelike, lightlike or timelike if $<x,x>_L > 0, <x,x>_L = 0$ or $<x,x>_L < 0$, respectively. The norm of the vector $x \in E_1^3$ is defined by

$$\left\| \sqrt{|<x,x>_L|} \right\|.$$

For any $x = (x_1, x_2, x_3), y = (y_1, y_2, y_3) \in E_1^3$, the Lorentzian vector product $x \times_L y$ is defined as follows:

$$x \times_L y = (-x_2 y_3 + x_3 y_2, x_3 y_1 - x_1 y_3, x_1 y_2 - x_2 y_1).$$

A curve $\alpha = \alpha(s): I \subset R \to E_1^3$ is spacelike, timelike or lightlike, if all of its tangent vectors $\alpha'(s)$ are respectively spacelike, timelike or lightlike for every $s \in I$.

Let $\alpha(s)$ be a unit speed non-lightlike curve in $E_1^3$ parametrized with the arclength s. Consider the Frenet frame $(T(s) = \alpha'(s), N(s), B(s))^T$ attached to the curve $\alpha(s)$ such that T is the unit tangent vector field, N is the principal normal vector field and B is the binormal vector field. The Frenet-Serret formulas are given by

$$\begin{bmatrix} T' \\ N' \\ B' \end{bmatrix} = \begin{bmatrix} 0 & \kappa & 0 \\ (\varepsilon_B)\kappa & 0 & \tau \\ 0 & (\varepsilon_T)\tau & 0 \end{bmatrix} \cdot \begin{bmatrix} T \\ N \\ B \end{bmatrix},$$

where $\varepsilon_X = <X,X>_L, \kappa = \kappa(s)$ and $\tau = \tau(s)$ are the curvature and torsion functions of the non-lightlike curve $\alpha(s)$, respectively (López, 2008; Choi et al. (2012).

## 3. Parallel frames of non-lightlike curves and vortex filament equation

The parallel frame is an alternative approach to defining a moving frame that is well-defined even when the curve has vanishing second derivative. We can parallel transport an orthonormal frame (such as Frenet frame) along a curve simply by parallel transporting each component of the frame. The parallel transport frame is based on the observation that, while the tangent vector for a given curve model is unique, we may choose any convenient arbitrary basis for the principal normal and binormal vector fields, so long as it is in the normal plane perpendicular to the tangent vector field at each point (Hanson and Ma, 1995). Here we describe parallel frames for timelike and spacelike curves in Minkowski 3-space.

### 3.1. The case of timelike curves in $E_1^3$

Let $\alpha(s): I \to E_1^3$ be a timelike unit speed curve. The Frenet frame of $\alpha$ can written as $(T, N, B)^T$ where $<T,T>_L = -1, <N,N>_L = 1, <B,B>_L = 1$ and $T \times N = B$. We want to change the Frenet frame $(T, N, B)^T$ to an orthonormal frame $(T, E_1, E_2)^T$ so that the 2,3-th and 3,2-th entries of the coefficient matrix of $(T, E_1, E_2)_s^T$ are zero. Throughout the paper, the subscripts s and t will denote differentiation with respect to these parameters, respectively. The following well known proposition is important for the construction of the parallel frames.

**Proposition 3.1.** Let $X_1, X_2, \ldots, X_n: [a,b] \to E_1^n$ be smooth maps and $(X_1, X_2, \ldots, X_n)$ be an orthonormal basis of $E_1^n$ for all $s \in [a,b]$. Then, there exists a matrix valued function $A(s) = \left(a_{ij}(s)\right)$ so that

(i) $\left(X_i(s)\right)_s = \sum_{j=1}^n a_{ij}(s) X_j(s)$,

(ii) $a_{ij} = \varepsilon_{X_i} < \left(X_j\right)_s, X_i >_L$,

(iii) $\varepsilon_{X_j} a_{ij} + \varepsilon_{X_i} a_{ji} = 0$ (that is $A(s)$ is semi-skew symmetric),

where $\varepsilon_{X_i} = <X_i, X_i>_L$.

Due to the method described in (Terng, 2005), rotating the frame $(N, B)$ by an angle $\theta(s)$ satisfying the equality $\left(\theta(s)\right)_s = \tau(s)$, we get the parallel frame $(E_1, E_2)$ with spacelike vector fields

$$\begin{cases} E_1 = cos\theta . N - sin\theta . B \\ E_2 = sin\theta . N + cos\theta . B \end{cases} \quad (1)$$

and $< (E_1)_s, E_2 >_L = < E_1, (E_2)_s >_L = 0$. So that the parallel frame $(T, E_1, E_2)^T$ satisfies

$$\begin{bmatrix} T \\ E_1 \\ E_2 \end{bmatrix}_s = \begin{bmatrix} 0 & k_1 & k_2 \\ k_1 & 0 & 0 \\ k_2 & 0 & 0 \end{bmatrix} . \begin{bmatrix} T \\ E_1 \\ E_2 \end{bmatrix} \quad (2)$$

where the principal curvatures of the curve $\alpha$ along $E_1$ and $E_2$ are respectively defined as

$$k_1 = \kappa cos\theta, \; k_2 = \kappa sin\theta. \quad (3)$$

However, the choice of the parallel frame is not unique because the angle $\theta(s)$ can be replaced by $\theta(s) + c, c \in R$. But still the frame is parallel and satisfies the equations in (2). Using (1) the relation between the Frenet frame and the parallel frame can be written as follows:

$$\begin{bmatrix} T \\ N \\ B \end{bmatrix} = \begin{bmatrix} 1 & 0 & 0 \\ 0 & cos\theta & sin\theta \\ 0 & -sin\theta & cos\theta \end{bmatrix} . \begin{bmatrix} T \\ E_1 \\ E_2 \end{bmatrix}. \quad (4)$$

The movement of a thin vortex in a viscous fluid can be modeled by the motion of a curve propagating in the direction of its binormal with a speed equal to its curvature according to the vortex filament equation

$$\alpha_t = \alpha_s \times \alpha_{ss}. \quad (5)$$

It is known that the binormal motions of spacelike and timelike curves in Minkowski 3-space are respectively equivalent to the non-linear heat system

$$\begin{cases} q_t = q_{ss} + q^2 r \\ r_t = -r_{ss} - r^2 q \end{cases} \quad (6)$$

where $r = \kappa e^{-\int_0^s \tau ds}, q = \kappa e^{\int_0^s \tau ds}$ and to the NLS$^-$ equation

$$iq_t + q_{ss} - \frac{1}{2} \|q\|^2 q = 0 \quad (7)$$

where $q = \kappa e^{i \int_0^s \tau ds}$ (Terng, 2005). The vortex filament equation has the following property (Hasimoto, 1972):

**Proposition 3.2.** If $\alpha(s, t)$ is a solution of the vortex filament equation and $\|\alpha_s(s, 0)\| = 1$ for all s, then $\|\alpha_s(s, t)\| = 1$ for all $(s, t)$. In other words, if $\alpha(\cdot, 0)$ is parametrized by arc length, then so is $\alpha(\cdot, t)$ for all t.

Now we can explain the geometric meaning of the evolution equation (5) on the space of timelike curves in $E_1^3$. Let $(T, N, B)^T(\cdot, t)$ denote the Frenet frame of the timelike curve $\alpha(\cdot, t)$. Then we get

$$\alpha_t = \alpha_s \times \alpha_{ss} = T \times T_s = T \times \kappa N = \kappa B.$$

That is, the motion of the timelike curve is in the direction of binormal with curvature as its speed. The expression of the equation (5) in terms of the parallel frame can be written using (3) and (4) as

$$\alpha_t = \kappa B = -k_2 E_1 + k_1 E_2.$$

In order to see the equivalence of the motion equation (5) and the NLS$^-$ equation (7), let $\alpha(s, t)$ be a solution of (5). We want to compute the time evolution of the parallel frame $(T, E_1, E_2)^T(\cdot, t)$. Using (2) and the previous equation we find

$$(T)_t = (\alpha_s)_t = (\alpha_t)_s = (-k_2 E_1 + k_1 E_2)_s$$
$$= -(k_2)_s E_1 - k_2(k_1 T) + (k_1)_s E_2 + k_1(k_2 T)$$
$$= -(k_2)_s E_1 + (k_1)_s E_2. \tag{8}$$

Since $(T, E_1, E_2)^T$ is orthonormal, according to Proposition 3.2 there exists a function $u(s, t)$ so that

$$\begin{bmatrix} T \\ E_1 \\ E_2 \end{bmatrix}_t = \begin{bmatrix} 0 & -(k_2)_s & (k_1)_s \\ -(k_2)_s & 0 & u \\ (k_1)_s & -u & 0 \end{bmatrix} \cdot \begin{bmatrix} T \\ E_1 \\ E_2 \end{bmatrix}$$

or equivalently

$$\begin{cases} (T)_t = -(k_2)_s E_1 + (k_1)_s E_2 \\ (E_1)_t = -(k_2)_s T + u E_2 \\ (E_2)_t = (k_1)_s T - u E_1. \end{cases} \tag{9}$$

Now to evaluate $u(s, t)$ we will use $u_s(s, t)$. With the equations (2,8,9) we obtain

$$< (E_1)_{st}, E_2 >_L = < (k_1 T)_t, E_2 >_L = k_1 < (T)_t, E_2 >_L$$
$$= k_1 < -(k_2)_s E_1 + (k_1)_s E_2, E_2 >_L$$
$$= k_1 (k_1)_s$$

$$< (E_1)_{ts}, E_2 >_L = < (-(k_2)_s T + u E_2)_s, E_2 >_L$$
$$= < -(k_2)_{ss} T - (k_2)_s (k_1 E_1 + k_2 E_2) + u_s E_2 + u(k_2 T), E_2 >_L$$
$$= -k_2 (k_2)_s + u_s.$$

So, by the equality of these equations we get

$$u_s = k_1(k_1)_s + k_2(k_2)_s$$

and we find

$$u(s,t) = \frac{1}{2}(k_1^2 + k_2^2) + c(t)$$

for some smooth function $c(t)$. Here, for each fixed t we can rotate $(E_1, E_2)(\cdot, t)$ by a constant angle $\theta(t)$ to another parallel frame $(\tilde{E}_1, \tilde{E}_2)(\cdot, t)$ of $\alpha(\cdot, t)$, i.e.,

$$\begin{cases} \tilde{E}_1(s,t) = cos\theta(t)E_1(s,t) + sin\theta(t)E_2(s,t), \\ \tilde{E}_2(s,t) = -sin\theta(t)E_1(s,t) + cos\theta(t)E_2(s,t). \end{cases}$$

Then by direct computation we can find that

$$\tilde{u} =< (\tilde{E}_1)_t, \tilde{E}_2 >_L = u + \theta_t = \frac{1}{2}(k_1^2 + k_2^2) + c(t) + \theta_t.$$

Taking $\theta_t = -c$ in the equation above, we have a new parallel frame which satisfies

$$\begin{bmatrix} T \\ \tilde{E}_1 \\ \tilde{E}_2 \end{bmatrix}_t = \begin{bmatrix} 0 & -(\tilde{k}_2)_s & (\tilde{k}_1)_s \\ -(\tilde{k}_2)_s & 0 & \frac{1}{2}(\tilde{k}_1^2 + \tilde{k}_2^2) \\ (\tilde{k}_1)_s & -\frac{1}{2}(\tilde{k}_1^2 + \tilde{k}_2^2) & 0 \end{bmatrix} \cdot \begin{bmatrix} T \\ \tilde{E}_1 \\ \tilde{E}_2 \end{bmatrix}.$$

Hence we have proved the first part of the following theorem:

**Theorem 3.3.** Let $\alpha(s,t)$ be a solution of the vortex filament equation (5) and $\|\alpha_s(s,0)\| = 1$ for all s. Then

(1) there exists a parallel normal frame $(T, E_1, E_2)^T(\cdot, t)$ for each timelike curve $\alpha(\cdot, t)$ so that

$$\begin{bmatrix} T \\ E_1 \\ E_2 \end{bmatrix}_s = \begin{bmatrix} 0 & k_1 & k_2 \\ k_1 & 0 & 0 \\ k_2 & 0 & 0 \end{bmatrix} \cdot \begin{bmatrix} T \\ E_1 \\ E_2 \end{bmatrix},$$

$$\begin{bmatrix} T \\ E_1 \\ E_2 \end{bmatrix}_t = \begin{bmatrix} 0 & -(k_2)_s & (k_1)_s \\ -(k_2)_s & 0 & \frac{1}{2}(k_1^2 + k_2^2) \\ (k_1)_s & -\frac{1}{2}(k_1^2 + k_2^2) & 0 \end{bmatrix} \cdot \begin{bmatrix} T \\ E_1 \\ E_2 \end{bmatrix}$$

where $k_1(\cdot, t)$ and $k_2(\cdot, t)$ are the principal curvatures of $\alpha(\cdot, t)$ along $E_1(\cdot, t)$ and $E_2(\cdot, t)$ respectively.

(2) $q = k_1 + ik_2$ is a solution of the NLS$^-$ equation

$$iq_t + q_{ss} - \frac{1}{2}\|q\|^2 q = 0.$$

**Proof.** First part is proved. For the second part we use the equations (2,9) to compute the time evolution equations of the principal curvatures $k_1$ and $k_2$:

$$(k_1)_t = (<T_s, E_1>_L)_t = <(T_s)_t, E_1>_L + <T_s, (E_1)_t>_L$$

$$= <(T_t)_s, E_1>_L + <T_s, (E_1)_t>_L = -(k_2)_{ss} + \frac{1}{2}(k_1^2 + k_2^2)k_2,$$

$$(k_2)_t = (<T_s, E_2>_L)_t = <(T_t)_s, E_2>_L + <T_s, (E_2)_t>_L$$

$$= (k_1)_{ss} - \frac{1}{2}(k_1^2 + k_2^2)k_1.$$

Substituting these in the NLS$^-$ equation we have the result:

$$iq_t + q_{ss} - \frac{1}{2}\|q\|^2 q = i((k_1)_t + i(k_2)_t) + ((k_1)_{ss} + i(k_2)_{ss}) - \frac{1}{2}(k_1^2 + k_2^2)(k_1 + ik_2) = 0.$$

*3.2. The case of spacelike curves with timelike principal normal in $E_1^3$*

Similar to the previous case we will explain the geometric meaning of the evolution equation (5) on the space of spacelike curves with timelike principal normal in $E_1^3$. Let $(T, N, B)^T(\cdot, t)$ denote the Frenet frame of the unit speed spacelike curve $\alpha(\cdot, t)$ so that $<T,T>_L = <B,B>_L = 1, <N,N>_L = -1$. If we rotate the frame $(N, B)$ with an angle $\theta(s)$ satisfying the equality $(\theta(s))_s = -\tau(s)$, we get the parallel frame $(E_1, E_2)$ with

$$\begin{cases} E_1 = cosh\theta . N + sinh\theta . B \\ E_2 = sinh\theta . N + cosh\theta . B \end{cases}$$

and $<(E_1)_s, E_2>_L = <E_1, (E_2)_s>_L = 0$. So that $E_1$ is timelike, $E_2$ is spacelike and the parallel frame $(T, E_1, E_2)^T$ satisfies

$$\begin{bmatrix} T \\ E_1 \\ E_2 \end{bmatrix}_s = \begin{bmatrix} 0 & k_1 & k_2 \\ k_1 & 0 & 0 \\ -k_2 & 0 & 0 \end{bmatrix} \cdot \begin{bmatrix} T \\ E_1 \\ E_2 \end{bmatrix} \tag{10}$$

where the principal curvatures of the curve $\alpha$ along $E_1$ and $E_2$ are respectively defined as

$$k_1 = \kappa cosh\theta, \quad k_2 = \kappa sinh\theta.$$

The relation between the Frenet frame and the parallel frame can be written as follows:

$$\begin{bmatrix} T \\ N \\ B \end{bmatrix} = \begin{bmatrix} 1 & 0 & 0 \\ 0 & cosh\theta & sinh\theta \\ 0 & sinh\theta & cosh\theta \end{bmatrix} \cdot \begin{bmatrix} T \\ E_1 \\ E_2 \end{bmatrix}.$$

Using (10) and the equation above, the vortex filament equation (5) can be expressed as

$$\alpha_t = \kappa B = k_2 E_1 + k_1 E_2.$$

In order to find the time evolution of the parallel frame $(T, E_1, E_2)^T$, first we find

$$(T)_t = (\alpha_s)_t = (\alpha_t)_s = (k_2 E_1 + k_1 E_2)_s = (k_2)_s E_1 + (k_1)_s E_2.$$

Again by the Proposition 3.2, since the frame $(T, E_1, E_2)^T$ is orthogonal there exists a function $u(s,t)$ so that

$$\begin{bmatrix} T \\ E_1 \\ E_2 \end{bmatrix}_t = \begin{bmatrix} 0 & (k_2)_s & (k_1)_s \\ (k_2)_s & 0 & u \\ -(k_1)_s & u & 0 \end{bmatrix} \cdot \begin{bmatrix} T \\ E_1 \\ E_2 \end{bmatrix}.$$

Similar to the case of timelike curves, to find the function $u(s,t)$ we proceed as follows:

$$< (E_1)_{st}, E_2 >_L = < (k_1 T)_t, E_2 >_L = k_1 < (T)_t, E_2 >_L$$

$$= k_1 < (k_2)_s E_1 + (k_1)_s E_2, E_2 >_L$$

$$= k_1 (k_1)_s$$

$$< (E_1)_{ts}, E_2 >_L = < ((k_2)_s T + u E_2)_s, E_2 >_L$$

$$= < (k_2)_{ss} T + (k_2)_s (k_1 E_1 + k_2 E_2) + u_s E_2 + u(-k_2 T), E_2 >_L$$

$$= k_2 (k_2)_s + u_s.$$

$$u_s = k_1 (k_1)_s - k_2 (k_2)_s \Rightarrow u(s,t) = \frac{1}{2}(k_1^2 - k_2^2) + c(t).$$

Hence, choosing the parallel frame such that $c(t)$ vanishes we can state the following theorem:

**Theorem 3.4.** Let $\alpha(s,t)$ be a solution of the vortex filament equation (5) and $\|\alpha_s(s,0)\| = 1$ for all s. Then

(1) there exists a parallel normal frame $(T, E_1, E_2)^T(\cdot, t)$ for each spacelike curve $\alpha(\cdot, t)$ with timelike principal normal so that

$$\begin{bmatrix} T \\ E_1 \\ E_2 \end{bmatrix}_s = \begin{bmatrix} 0 & k_1 & k_2 \\ k_1 & 0 & 0 \\ -k_2 & 0 & 0 \end{bmatrix} \cdot \begin{bmatrix} T \\ E_1 \\ E_2 \end{bmatrix},$$

$$\begin{bmatrix} T \\ E_1 \\ E_2 \end{bmatrix}_t = \begin{bmatrix} 0 & (k_2)_s & (k_1)_s \\ (k_2)_s & 0 & \frac{1}{2}(k_1^2 - k_2^2) \\ -(k_1)_s & \frac{1}{2}(k_1^2 - k_2^2) & 0 \end{bmatrix} \cdot \begin{bmatrix} T \\ E_1 \\ E_2 \end{bmatrix}$$

where $k_1(\cdot, t)$ and $k_2(\cdot, t)$ are the principal curvatures of $\alpha(\cdot, t)$ along $E_1(\cdot, t)$ and $E_2(\cdot, t)$ respectively.

(2) $\begin{cases} q = \frac{1}{\sqrt{2}}(k_1 + k_2) \\ r = -\frac{1}{\sqrt{2}}(k_1 - k_2) \end{cases}$ is a solution of the non-linear heat system

$$\begin{cases} q_t = q_{ss} + q^2 r \\ r_t = -r_{ss} - r^2 q. \end{cases}$$

**Proof.** For the second part, first we compute the time evolution equations of the principal curvatures $k_1$ and $k_2$:

$$(k_1)_t = (-<T_s, E_1>_L)_t = -<(T_t)_s, E_1>_L -<T_s, (E_1)_t>_L$$

$$= (k_2)_{ss} - \frac{1}{2}(k_1^2 - k_2^2)k_2,$$

$$(k_2)_t = (<T_s, E_2>_L)_t = <(T_t)_s, E_1>_L +<T_s, (E_1)_t>_L$$

$$= (k_1)_{ss} - \frac{1}{2}(k_1^2 - k_2^2)k_1.$$

Then substituting these equations in the non-linear heat system we finish the proof:

$$q_t = \frac{1}{\sqrt{2}}((k_1)_t + (k_2)_t) = \frac{1}{\sqrt{2}}\left((k_2)_{ss} - \frac{1}{2}(k_1^2 - k_2^2)k_2 + (k_1)_{ss} - \frac{1}{2}(k_1^2 - k_2^2)k_1\right) = q_{ss} + q^2 r$$

$$r_t = -\frac{1}{\sqrt{2}}((k_1)_t - (k_2)_t) - \frac{1}{\sqrt{2}}\left((k_2)_{ss} - \frac{1}{2}(k_1^2 - k_2^2)k_2 - (k_1)_{ss} + \frac{1}{2}(k_1^2 - k_2^2)k_1\right) = -r_{ss} - r^2 q.$$

*3.3. The case of spacelike curves with timelike binormal in $E_1^3$*

Finally we will describe the geometric meaning of the evolution equation (5) on the space of spacelike curves with timelike binormal in $E_1^3$. Let $(T, N, B)^T(\cdot, t)$ denote the Frenet frame of the unit speed spacelike curve $\alpha(\cdot, t)$ so that $<T, T>_L = <N, N>_L = 1, <B, B>_L = -1$. Rotating the frame $(N, B)$ with an angle $\theta(s)$ such that $(\theta(s))_s = -\tau(s)$ is satisfied, we will get the parallel frame $(E_1, E_2)$. Here $E_1$ is spacelike, $E_2$ is timelike and the parallel frame $(T, E_1, E_2)^T$ satisfies

$$\begin{bmatrix} T \\ E_1 \\ E_2 \end{bmatrix}_s = \begin{bmatrix} 0 & k_1 & k_2 \\ -k_1 & 0 & 0 \\ k_2 & 0 & 0 \end{bmatrix} \cdot \begin{bmatrix} T \\ E_1 \\ E_2 \end{bmatrix}$$

where the principal curvatures of the curve $\alpha$ along $E_1$ and $E_2$ are respectively defined as $k_1 = \kappa \cosh\theta$, $k_2 = \kappa \sinh\theta$. The expression of the equation (5) and the time evolution of the tangent vector field T are same as in the previous case. By the Proposition 3.2 there exists a function $u(s, t)$ so that

$$\begin{bmatrix} T \\ E_1 \\ E_2 \end{bmatrix}_t = \begin{bmatrix} 0 & (k_2)_s & (k_1)_s \\ -(k_2)_s & 0 & u \\ (k_1)_s & u & 0 \end{bmatrix} \cdot \begin{bmatrix} T \\ E_1 \\ E_2 \end{bmatrix}.$$

Then we can find $u(s,t)$ as follows:

$$< (E_1)_{st}, E_2 >_L = < (-k_1 T)_t, E_2 >_L = -k_1 < (T)_t, E_2 >_L = k_1 (k_1)_s$$

$$< (E_1)_{ts}, E_2 >_L = < ((k_2)_s T + u E_2)_s, E_2 >_L = k_2 (k_2)_s - u_s.$$

$$u_s = -k_1(k_1)_s + k_2(k_2)_s \Rightarrow u(s,t) = -\frac{1}{2}(k_1^2 - k_2^2) + c(t).$$

Hence, we can state the following theorem for this final case:

**Theorem 3.5.** Let $\alpha(s,t)$ be a solution of the vortex filament equation (5) and $\|\alpha_s(s,0)\| = 1$ for all s. Then

(1) there exists a parallel normal frame $(T, E_1, E_2)^T(\cdot, t)$ for each spacelike curve $\alpha(\cdot, t)$ with timelike binormal so that

$$\begin{bmatrix} T \\ E_1 \\ E_2 \end{bmatrix}_s = \begin{bmatrix} 0 & k_1 & k_2 \\ -k_1 & 0 & 0 \\ k_2 & 0 & 0 \end{bmatrix} \cdot \begin{bmatrix} T \\ E_1 \\ E_2 \end{bmatrix},$$

$$\begin{bmatrix} T \\ E_1 \\ E_2 \end{bmatrix}_t = \begin{bmatrix} 0 & (k_2)_s & (k_1)_s \\ -(k_2)_s & 0 & -\frac{1}{2}(k_1^2 - k_2^2) \\ (k_1)_s & -\frac{1}{2}(k_1^2 - k_2^2) & 0 \end{bmatrix} \cdot \begin{bmatrix} T \\ E_1 \\ E_2 \end{bmatrix}$$

where $k_1(\cdot, t)$ and $k_2(\cdot, t)$ are the principal curvatures of $\alpha(\cdot, t)$ along $E_1(\cdot, t)$ and $E_2(\cdot, t)$ respectively.

(2) $\begin{cases} q = \frac{1}{\sqrt{2}}(k_1 + k_2) \\ r = \frac{1}{\sqrt{2}}(k_1 - k_2) \end{cases}$ is a solution of the non-linear heat system

$$\begin{cases} q_t = q_{ss} + q^2 r \\ r_t = -r_{ss} - r^2 q. \end{cases}$$

**Proof.** First of all we compute the time evolution equations of the principal curvatures $k_1$ and $k_2$:

$$(k_1)_t = (< T_s, E_1 >_L)_t = (k_2)_{ss} + \frac{1}{2}(k_1^2 - k_2^2)k_2,$$

$$(k_2)_t = (-< T_s, E_2 >_L)_t = (k_1)_{ss} + \frac{1}{2}(k_1^2 - k_2^2)k_1.$$

Then substituting these in the non-linear heat system we can easily find that the given is a solution of this system.

**References**

Choi, J.H., Kimb, Y.H. and Ali, A.T. (2012). Some associated curves of Frenet non-lightlike curves in $E_1^3$. *J. Math. Anal. Appl.*, 394, 712-723.

Ding, Q. (1998). A note on the NLS and the Schrödinger flow of maps. *Physics Letters A*, 248, 49-56.

Ding, Q. and Inoguchi, J.-i. (2004). Schrödinger flows, binormal motion for curves and the second AKNS-hierarchies. *Chaos, Solitons & Fractals*, 21, 669-677.

Gürbüz, N. (2010). Intrinsic Geometry of the NLS Equation and Heat System in 3-Dimensional Minkowski Space. Adv. Studies Theor. Phys., 4(11), 557-564.

Hanson, A.J. and Ma, H. (1995). Technical Report TR425: Parallel Transport Approach to Curve Framing, Department of Computer Science, Indiana University.

Hasimoto, H. (1972). A solution on a vortex filament. J. Fluid Mech., 51, 477-485.

López, R. (2008). Differential Geometry of Curves and Surfaces in Lorentz-Minkowski space. arXiv:0810.3351 [math.DG]

Ricca, R. L. (1996). The contributions of Da Rios and Levi-Civita to asymptotic potential theory and vortex filament dynamics. Fluid Dynamics Research, 18, 245-268.

Rogers, C. and Schief, W. K. (1998). Intrinsic Geometry of the NLS Equation and Its Auto-Bäcklund Transformation. *Studies in Applied Mathematics*, 101, 267-287.

Terng, C. (2005). Lecture 1 Peking University Summer School.

Zhang H. & Wu, F. (2007). Vortex Filament Equation and Non-linear Schrödinger Equation in S3. Kyungpook Math. J., 47, 381-392.